\documentclass[9pt,bezier]{article}
\usepackage{amsmath,amssymb,amsfonts,euscript}
\usepackage{algorithmic, algorithm}
\usepackage{graphicx}
\usepackage{xcolor}

\usepackage[ansinew]{inputenc}
\usepackage{graphicx}
\usepackage{color}
\usepackage{mathrsfs}

\usepackage[colorlinks]{hyperref}
\usepackage[active,new,noold,marker]{xrcs}
\catcode`\=13
\def{$\bowtie$}

\usepackage{eurosym}
\usepackage{lscape} 

\newtheorem{prethm}{{\bf Theorem}}

\newenvironment{thm}{\begin{prethm}{\hspace{-0.5
               em}{\bf}}}{\end{prethm}}

\newtheorem{prepro}{{\bf Theorem}}

\newtheorem{preprop}{{\bf Proposition}}

\newtheorem{precor}{{\bf Corollary}}

\newtheorem{preconj}{{\bf Conjecture}}

\newtheorem{predefi}{{\bf Definition}}

\newtheorem{preremark}{{\bf Remark}}

\newenvironment{remark}{\begin{preremark}\rm{\hspace{-0.5
               em}{\bf}}}{\end{preremark}}

\newtheorem{preexample}{{\bf Example}}

\newtheorem{prelem}{{\bf Lemma}}

\newtheorem{prelam}{{\bf Lemma}}

\newtheorem{preprob}{{\bf Problem}}

\newenvironment{prob}{\begin{preprob}{\hspace{-0.5
               em}{\bf.}}}{\end{preprob}}

\newtheorem{preproof}{{\bf Proof}}

\newtheorem{preali}{{\bf Proof of Theorem 1.}}

\newenvironment{ali}[1]{\begin{preali}{\rm
               #1}\hfill{$\Box$}}{\end{preali}}

\newtheorem{prealii}{{\bf Proof of Theorem 2.}}

\newenvironment{alii}[1]{\begin{prealii}{\rm
               #1}\hfill{$\Box$}}{\end{prealii}}

\newtheorem{prealiii}{{\bf Proof of Theorem 3.}}

\newenvironment{aliii}[1]{\begin{prealiii}{\rm
               #1}\hfill{$\Box$}}{\end{prealiii}}

\newtheorem{prealiiii}{{\bf Proof of Theorem 4.}}

\newenvironment{aliiii}[1]{\begin{prealiiii}{\rm
               #1}\hfill{$\Box$}}{\end{prealiiii}}

\newtheorem{prealij}{{\bf Proof of Theorem 5.}}

\newtheorem{prealijj}{{\bf Proof of Theorem 6.}}

\newtheorem{prealijjj}{{\bf Proof of Theorem 7.}}

\newtheorem{prealijjjk}{{\bf Proof of Theorem 8.}}

\title{Is there any polynomial upper bound
for the universal labeling of graphs?}

\author{{\normalsize
{  Arash Ahadi${}^{\mathsf{a}}$},\,
 {  Ali Dehghan${}^{\mathsf{b}}$},\,
{  Morteza Saghafian${}^{\mathsf{a}}$},\,
}\vspace{3mm}
\\{\footnotesize{${}^{\mathsf{a}}$\it Department of
Mathematical Sciences, Sharif University of Technology, Tehran,
Iran}}  {\footnotesize{}}\\{\footnotesize{${}^{\mathsf{b}}$\it
Systems and Computer Engineering Department, Carleton University, Ottawa,   Canada}}
\thanks{{\it E-mail addresses}:  $\mathsf{arash\_ahadi@mehr.sharif.edu}$ (Arash Ahadi),  $\mathsf{alidehghan@sce.carleton.ca}$ (Ali Dehghan),
$\mathsf{saghafian@ce.sharif.edu}$ (Morteza Saghafian). } }

\begin{document}
\maketitle

\begin{abstract}
{\small \noindent
A {\it universal labeling} of a graph $G$ is a labeling of the edge set in $G$ such that in every orientation $\ell$ of $G$ for every two adjacent vertices $v$ and $u$, the sum of incoming edges of $v$ and $u$ in the oriented graph are different from each other. The {\it universal labeling number} of a graph $G$ is the minimum number $k$ such that $G$ has {\it universal labeling} from $\{1,2,\ldots, k\}$ denoted it by $\overrightarrow{\chi_{u}}(G) $. We have
$2\Delta(G)-2 \leq \overrightarrow{\chi_{u}} (G)\leq  2^{\Delta(G)}$, where $\Delta(G)$ denotes  the maximum
degree of $G$. In this work, we offer a provocative question that is:" Is there  any polynomial function $f$ such that for every graph $G$, $\overrightarrow{\chi_{u}} (G)\leq f(\Delta(G))$?".
Towards this question, we introduce some lower and upper bounds on their parameter of
interest. Also, we prove that for every tree $T$, $\overrightarrow{\chi_{u}}(T)=\mathcal{O}(\Delta^3) $.
Next, we show that for a given 3-regular graph $G$, the  universal labeling number of $G$ is 4 if and only if $G$ belongs to Class 1.
Therefore,  for a given 3-regular graph $G$, it is an $ \mathbf{NP} $-complete  to determine whether
 the  universal labeling number of $G$ is 4.
Finally, using probabilistic methods, we almost confirm a weaker version of the problem.

}

\begin{flushleft}

\end{flushleft}

\end{abstract}


\section{Introduction}
\label{}

Throughout the paper we denote $\{1,2,\ldots, k\}$ by $\mathbb{N}_k$
and we use \cite{MR1367739} for terminology
and notations which are not defined here, also we consider only simple and finite graphs and digraphs.
In 2004 Karo\'nski, \L{}uczak and Thomason
initiated the study of
{\it neighbour-sum-distinguishing} labeling \cite{MR2047539}. They introduced an edge-labeling which
is additive vertex-coloring that means for every edge
$uv$, the sum of labels of the edges incident to $u$ is different
from the sum of labels of the edges incident to $v$.
It was conjectured in  \cite{MR2047539} that every graph with no isolated edge has a {\it neighbour-sum-distinguishing} labeling
from $\mathbb{N}_3$ (1-2-3-conjecture). This conjecture has been studied extensively by several
authors, for instance see \cite{iapprox, MR3072733,  MR2047539}.
Currently, we know that  every connected graph with more than two vertices has a {\it neighbour-sum-distinguishing} labeling, using the labels from $\mathbb{N}_{5}$ \cite{MR2595676}.

Various directed  versions of the problem of vertex distinguishing colorings are of interest
and have been recently investigated more intensively, see for instance \cite{MR3095464, MR3293286, MR3514378}.
In this work we consider a new directed version of this problem.
A {\it universal labeling} of a graph $G$ is a labeling of the edge set in $G$ such that in every orientation $\ell$ of $G$ for every two adjacent vertices $v$ and $u$, the sum of incoming edges of $v$ and $u$ in the oriented graph are different from each other. The {\it universal labeling number} of a graph $G$ is the minimum number $k$ such that $G$ has {\it universal labeling} from $\mathbb{N}_k$, denoted it by $\overrightarrow{\chi_{u}}(G) $.
For this notion, we conjecture that
every graph should be weightable with a polynomial number of weights
only, where the notion of polynomiality is with respect to the graph
maximum degree parameter. As first steps towards the conjecture,
we introduce some lower and upper bounds on their parameter of
interest, then we prove the conjecture for trees, and partially for 3-regular
and 4-regular graphs, and finally using probabilistic methods, we almost confirm a weaker version of the conjecture.

\section{Universal Labeling}

\begin{figure}[ht]
\begin{center}
\includegraphics[scale=.6]{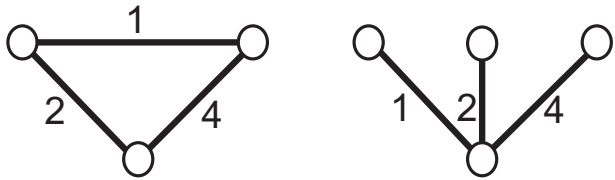}
\caption{A universal labeling for the graph $G=K_3 \cup K_{1,3}$.} \label{graphA}
\end{center}
\end{figure}

A {\it universal labeling} of a graph $G$ is a labeling of the edge set in $G$ such that in every orientation $\ell$ of $G$ for every two adjacent vertices $v$ and $u$, the sum of incoming edges of $v$ and $u$ in the oriented graph are different from each other. The {\it universal labeling number} of a graph $G$ is the minimum number $k$ such that $G$ has {\it universal labeling} from $\mathbb{N}_k$, denoted it by $\overrightarrow{\chi_{u}}(G) $. See Figure \ref{graphA}

Every graph has some {\it universal labelings}, for example one may put the different powers of
two $\{2^{i-1}:   i \in \mathbb{N}_{n}\}$ on the edges of $G$, where $n$ is the number of vertices.
Motivated by 1-2-3-conjecture, the following question posed by the second author on MathOverflow \cite{math33}.

\begin{prob}\label{p1}
Is there  a polynomial function $f$ such that for every graph $G$, $\overrightarrow{\chi_{u}} (G)\leq f(\Delta(G))$?
\end{prob}

Let $G$ be a graph and $e$ be an arbitrary edge in $G$. Without lose of generality suppose that $\ell$ is a universal labeling for $G$, it is easy to see that $\ell$ is also a universal labeling for $G\setminus \{e\}$.
Thus, in the set of graphs with $n$ vertices, complete graph $K_n$ has the maximum universal labeling number, so we have the following problem:
Is there  a polynomial function $f$ in terms of $n$ such that $\overrightarrow{\chi_{u}} (K_n)\leq f(n)$?
If the answer to this problem is no, then there is a graph $G$ and an edge $e\in E(G)$ such that the universal labeling number of $G\setminus \{e\}$ is a polynomial  in terms of $n$, but the universal labeling number of $G$ is not  polynomial  in terms of $n$. Finding such graphs can be interesting.

For $ k\in \mathbb{N} $, a {\it proper edge $k$-coloring} of $G$ is a function $c:
E(G)\rightarrow \mathbb{N}_k$, such that if $e,e'\in E(G)$ share a common endpoint,
then $c(e)$ and $c(e')$ are different.
The smallest integer $k$ such that the graph
$G$ has a proper edge $k$-coloring is called the {\it chromatic index} of the graph $G$ and denoted by $\chi '(G)$.
Let $f $ be a proper edge coloring for a given graph $G$. Then the function $\ell(e)= 2^{f(e)-1}$ is a universal labeling for the graph $G$. By Vizing's theorem \cite{MR0180505}, the chromatic index of a graph $G$ is equal to either $ \Delta(G) $ or $ \Delta(G) +1 $. So, every graph $G$ has a universal labeling from $\{2^{i-1}:i\in \mathbb{N}_{\Delta(G)+1}\}$. On the other hand,  note that every universal labeling for the edges of $G$ is a  proper edge coloring of $G$. Therefore the universal labeling number is at least the chromatic index of a graph. Thus,
$$\Delta(G) \leq \overrightarrow{\chi_{u}} (G)\leq  2^{\Delta(G)}.$$

Let $G$ be a graph and $f: E(G) \rightarrow  \mathbb{N}_{\overrightarrow{\chi_{u}} (G)}$ be a universal labeling for it.
Suppose that $v$ is a vertex with maximum degree $\Delta$ and let $\{ v_1,\ldots, v_{\Delta}\}$ be the set of neighbors of the vertex $v$. Since every universal labeling for the edges of $G$ is a  proper edge coloring, thus  $f(vv_1), f(vv_2), \ldots ,f(vv_{\Delta})$ are distinct numbers. Without loss of generality suppose that $f(vv_1) < f(vv_2)< \ldots <f(vv_{\Delta})=M$. Now, consider the following partition for $\mathbb{N}_{M-1}$:
$\big\{\{i,M-i\}: i \in \mathbb{N}_{\lfloor M/2 \rfloor} \big\} $.
The set $\{f(vv_i):i \in \mathbb{N}_{\Delta-1} \}$ contains at most one number from each of the above sets, so
$\Delta(G) \leq \lfloor M/2 \rfloor +1  $, therefore,

 \begin{equation}
2\Delta(G)-2 \leq \overrightarrow{\chi_{u}} (G)\leq  2^{\Delta(G)}.
\end{equation}

\begin{remark} Here we show that the lower bound in Eq. 1 is sharp.
Consider the complete bipartite graph $K_{1,m}$, by Eq. 1, $2m-2 \leq \overrightarrow{\chi_{u}} (G)$. Consider an arbitrary order for the  edges of $G$ and label them by $m-1,m,\ldots,2m-2$, respectively. This labeling is a universal labeling. Thus $  \overrightarrow{\chi_{u}} (K_{1,m})=2m-2$.
\end{remark}

In the first theorem, we solve Problem \ref{p1}, for trees.

\begin{thm}\label{pp1}
For every tree $T$, $\overrightarrow{\chi_{u}}(T)=\mathcal{O}(\Delta^3) $, where $\Delta$ is the maximum degree of the graph.
\end{thm}

Let $G$ be a connected graph and $\overrightarrow{\chi_{u}} (G) \leq 3$. By Eq. 1, $\Delta(G)\leq 2$. Therefore $G$ is a cycle or path (note that there is no graph with $\overrightarrow{\chi_{u}} = 3$). For $\overrightarrow{\chi_{u}} (G) \leq 4$, we have the following hardness.

\begin{thm}\label{pp2} Let $G$ be a 3-regular graph, then   $\overrightarrow{\chi_{u}} (G)= 4$ if and only if $G$ belongs to Class 1. Also, for a given 4-regular graph $G$, $\overrightarrow{\chi_{u}} (G)= 7$ if and only if $G$ belongs to Class 1.
\end{thm}

It was shown in  \cite{MR689264} that it is $ \mathbf{NP} $-hard to determine the chromatic index of an
$r$-regular graph for any $r \geq 3$, therefore,  for a given 3-regular graph $G$, it is an $ \mathbf{NP} $-complete  to determine whether the  universal labeling number of $G$ is 4. Also, for a given 4-regular graph $G$, it is an $ \mathbf{NP} $-complete  to determine whether the  universal labeling number of $G$ is 7.

An edge labeling $\ell$ for a graph $G$ is {\it almost every where universal labeling} if for a random orientation $\overrightarrow{O}$, $\Pr(\ell \text{ over }\overrightarrow{O}\text{ is proper})=1$. In other words, for the labeling $\ell$ of $G$ for every two adjacent vertices $v$ and $u$, the sum of incoming edges of $v$ and $u$ in the oriented graph which is obtained from the random orientation $\overrightarrow{O}$ are different from each other. Although we do not know any polynomial upper bound  for the universal labeling number, but there is a Quasi-polynomial upper bound for  almost every where universal labeling.

\begin{thm}
Every graph $G$ has an almost every where universal labeling from $\mathbb{N}_{n^{\frac{\lg n}{\lg \lg n}}}$, where $n$ is the number of vertices of the graph.
\end{thm}

\begin{remark}
For the complete graph $K_n$, let $f(n)=\omega(n^2)$.
Label every edge of $K_n$  randomly and independently by one color from $\mathbb{N}_{f(n)}$ with the
same probability.
Let $\overrightarrow{O}$ be a random orientation for the graph. The probability that the two ends of an arbitrary edges $e$, have equal sums of incoming edges (mod f(n)), is $1/f(n)$ . So the probability that this labeling is proper over $\overrightarrow{O}$ is $1- \mathcal{O}(\frac{n^2}{f(n)})\cong 1$.
Note that although the universal labeling number  is an increasing property, but one cannot consider the above computing to obtain an upper bound for almost every where universal labeling number of all graphs.
\end{remark}

{\bf The universal labeling game}

{\it A universal labeling game} is a perfect-information game played between two
players. The input of the game is an undirected graph  $G$ and number $k$ (the input is denoted by $(G,k)$). The game consists of $|E(G)|$ rounds. In each round
the first player chooses an unlabeled edge $e$ and label it with one of the numbers $\mathbb{N}_k$, afterwards, the next player chooses an orientation for $e$. Let $f$ be the labeling of $G$ after  $|E(G)|$ rounds and $D$ be its oriented graph.
If every two adjacent vertices of $D$ receive
distinct sums of  incoming labels in $f$, then first player wins.
{\it The universal labeling game number} of $G$, denoted by $\overrightarrow{\chi_{u}^{g}} (G) $, is the minimum number $k$ such that the first player has a wining strategy on $(G,k)$. It is clear that $\overrightarrow{\chi_{u}^{g}} (G) \leq \overrightarrow{\chi_{u}} (G) $, we prove the following upper bounds for $G$.

\begin{thm}\label{pp4} \\
{\em (i)} For a connected graph $G$, $ \overrightarrow{\chi_{u}^{g}} (G)=2$ if and if $G$ is a path or even cycle.\\
{\em (ii)} For every graph $G$, $ \overrightarrow{\chi_{u}^{g}} (G) \leq 2\Delta(G)$.\\
{\em (iii)} For every tree $T$, $\Delta(T) \leq  \overrightarrow{\chi_{u}^{g}} (T) \leq \Delta(T)+1$.
\end{thm}

In the above theorem we show that for every tree $T$, $\Delta(T) \leq  \overrightarrow{\chi_{u}^{g}} (T) \leq \Delta(T)+1$.
For a given tree $T$, determining the complexity of computing $\overrightarrow{\chi_{u}^{g}} (T)$  is interesting.
Finally, we ask the following about the computational complexity of the universal labeling game number of a given graph $G$.

\begin{prob}{
Is computing of the universal labeling game number of a given graph $\mathbf{NP} $-h?
}\end{prob}

\section{Proofs}

\begin{ali}{
Let $T$ be a tree with maximum degree $\Delta$ and $c:E(G)\rightarrow \mathbb{N}_{\Delta}$ be a proper edge coloring for $T$
(note that the chromatic index of every tree $T$ is $\Delta$ \cite{MR1367739}).
Define the function $N(e)=c(e)+\Delta-2$. We say that a set of numbers is sum-free if for every $i>1$, the sum of each $i$ members of that set is not member of it. Note that for each vertex $v$ the set of numbers $\{N(e): e\ni v\}$ is sum-free ({\bf Fact 1}).

Choose an arbitrary vertex $z$ of $T$, and perform a breadth-first search
algorithm from the vertex $z$. This defines a partition $V_0, V_1, \ldots, V_d$ of the vertices of $T$
where each part $V_i$ contains the vertices of $T$ which are at depth $i$ (at distance exactly $i$
from $z$). Each edge is between two parts $i-1$ and $i$ for a natural number $i$. Say an edge is at even-level if $i-1$ is even and odd-level otherwise. Each vertex at even depth is incident with just one odd-level edge  (except the root which is not incident with any odd-level edge) and each vertex at odd depth is incident with just one even-level edge ({\bf Fact 2}).

Let $K=(\Delta-1)+\Delta+(\Delta+1)+\cdots+(2\Delta-2)=\frac{3}{2}\Delta(\Delta-1) $. Now put the labels on the edges as follows: If an edge $e$ is at even-level, then put $N(e)$ on it and	if an edge $e$ is at odd-level, then put $N(e)\times K$ on it.
We claim that, this is a universal labeling from $\{\Delta-1,\Delta,\ldots,3\Delta(\Delta-1)^2\}$.

For every orientation of $T$ and for each vertex $v$, if we divide the sum of incoming edges of $v$ to $K$, we have $K\times c(v)+r(v)$, where $c(v)$ is the sum of $N(e)$'s for odd-level incoming edges of $v$ and $r(v)$ is the sum of $N(e)$'s for even-level incoming edges of $v$ ({\bf Fact 3}).

Consider an orientation on $T$. Suppose an edge $e$ with vertices $v_{i-1}$, $v_i$ at depth $i-1$, $i$  respectively, such that the sum of labels of incoming edges of $v_{i-1}$ and $v_i$ are equal. We have two cases:

Case 1. $i-1$ is even.
By Fact 3 and because $e$ is an incoming edge of $v_{i-1}$ or $v_i$, we have $r(v_i )=r(v_{i-1} )>0$. But by Fact 2, $v_i$ has only $e$ as an even-level edge. Therefore $N(e)$ should be equal to sum of $N(e')$'s for incoming even-level edges of $v_{i-1}$. By Fact 1, this is a contradiction.

Case 2. $i-1$ is odd.
By Fact 3 and since $e$ is an incoming edge of $v_{i-1}$ or $v_i$ we have $c(v_i )=c(v_{i-1} )>0$. But by Fact 2, $v_i$ has only $e$ as an odd-level edge. Therefore $N(e)$ should be equal to sum of $N(e')$'s for incoming odd-level edges of $v_{i-1}$, by Fact 1, this is a contradiction.

Thus that labeling is a universal labeling with $\mathcal{O}(\Delta^3)$ labels. This completes the proof.
}\end{ali}

\begin{alii}{

Let $G$ be a given 3-regular graph and $\ell: E(G)\rightarrow \mathbb{N}_4$ be a universal labeling for it. Also, suppose that $v$ is an arbitrary vertex and $u_1,u_2,u_3$ are its neighbors. First,
note that every universal labeling for the edges of $G$ is a  proper edge coloring. Therefore the universal labeling number is at least the chromatic index of a graph. So, $|\{\ell(vu_i):i\in \mathbb{N}_3\}|=3$. Also, note that the set of numbers $\{\ell(vu_i):i\in \mathbb{N}_3\}$ is sum-free, thus this set is not $123$ or $134$. Therefore, the set of numbers $\{\ell(vu_i):i\in \mathbb{N}_3\}$ is $124$ or $234$.

First, suppose that $G$ has the chromatic index 3 and let $f $ be a such proper edge coloring with colors $0,1,2$. The function
$2^f$ is a universal labeling for its edges. Next, assume that $G$ has a universal labeling $\ell$.
For every vertex $v$ and its neighbors $u_1,u_2,u_3$, the set of numbers $\{\ell(vu_1),\ell(vu_2),\ell(vu_3)\}$ is $124$ or $234$. So the following function is a proper edge coloring for $G$.

$ f(v)=
\begin{cases}
   \ell(v),      & $   if  $\,\, \ell(v)\neq 3,\\
   1,            & $   if  $\,\, \ell(v)= 3.
\end{cases}$
\\
\\
This completes the proof for the first part of theorem.

Now, assume that $G$ is a 4-regular graph. by Eq. 1, $ \overrightarrow{\chi_{u}} (G)\geq 6$.
First, we show that $ \overrightarrow{\chi_{u}} (G)\geq 7$. To the contrary, suppose that $ \overrightarrow{\chi_{u}} (G)= 6$ and let $f:E(G) \rightarrow \mathbb{N}_6$ be a universal labeling for $G$. Suppose that there is vertex $v\in V(G)$ such that $6 \notin \{f(uv)| uv \in E(G)\}$. So, $\{f(uv)| uv \in E(G)\} \subset \mathbb{N}_5$. It is clear that $1\notin \{f(uv)| uv \in E(G)\}$, therefore, the set of numbers $\{f(uv)| uv \in E(G)\}$ is $2345$. But that set of numbers is not sum-free, so this is a contradiction.
Consequently, for every vertex $v\in V(G)$, $6 \in \{f(uv)| uv \in E(G)\}$.
Since $6 \in \{f(uv)| uv \in E(G)\}$, thus, at most, one of the numbers 2,4 is in $\{f(uv)| uv \in E(G)\}$, Also, one of the numbers 1,5 is in $\{f(uv)| uv \in E(G)\}$. Therefore, $3 \in \{f(uv)| uv \in E(G)\}$. There are four cases for the set $ \{f(uv)| uv \in E(G)\}$: $6321,6341,6325$ and $6345$. The three former cases are not sum-free, so the set of numbers for every vertex $v$ is $ \{f(uv)| uv \in E(G)\}$ is $6345$.
Let $vu\in E(G)$, Since $3+6=4+5$ therefore, the two adjacent vertices $v$ and $u$ can have equal indegree, that is a contradiction.
So, $ \overrightarrow{\chi_{u}} (G)\geq 7$. Next, we show that every vertex $v$, $7\in \{f(vw)| vw \in E(G)\}$. Suppose that there is vertex $v$ such that $7\notin \{f(wv)| wv \in E(G)\}$. Without loss of generality suppose that $vu\in E(G)$ and $f(uv)=4$. We have $7\in \{f(wu)| wu \in E(G)\}$. In this situation $7=4+3$ and this a contradiction.
Therefore, every vertex $v$, $7\in \{f(uv)| uv \in E(G)\}$.
One of the numbers 1,6 is in $\{f(uv)| uv \in E(G)\}$, also, one of the numbers 2,5 is in $\{f(uv)| uv \in E(G)\}$ and one of the numbers 3,4 is in $\{f(uv)| uv \in E(G)\}$. Hence, there are eight cases for $\{f(uv)| uv \in E(G)\}$. Among those sets only $1357,4567,3567,2367 $ are sum-free.

Since for every vertex $v$, $7\in \{f(uv)| uv \in E(G)\}$, thus if $a+b=c+7$, then the set of numbers $\{f(uv)| uv \in E(G)\}$ is not $abc7$. By this fact the set of numbers $\{f(uv)| uv \in E(G)\}$ is not $1357,4567,2367 $. Thus $\{f(uv)| uv \in E(G)\}$
is $3567$. So $f$ is a proper 4-edge coloring for $G$. On the other hand if $G$ belongs to class 1, then let $\ell $ be one of its proper 4-edge coloring with colors $3567$. This coloring is also a universal labeling for $G$. This completes the proof.
}\end{alii}

Next, we prove that every graph $G$ has an almost every where universal labeling from $\mathbb{N}_{n^{\frac{\lg n}{\lg \lg n}}}$, where $n$ is the number of vertices of the graph. The key idea is: partitioning the edges of the graph into two parts based on the degrees of the vertices and labeling each part independently.

\begin{aliii}{
Let $G$ be a graph with $n$ vertices. Without loss of generality assume that  $\frac{\lg n }{\lg \lg n}$ is an integer number (if $\frac{\lg n }{\lg \lg n}$ is not  integer, we can consider $\lfloor \frac{\lg n }{\lg \lg n}\rfloor$ in the proof).
Consider the edge-induced subgraph of $G$ over all edges like $e$ such that both ends of $e$ have degrees less than $\lg n$ in $G$ and call it $H$.
Let $f$ be a proper edge coloring  of $H$ by labels $\{2^{i-1}: i \in \mathbb{N}_{\lg n}\}$.
Now, assume that $ 2n < p_1 < p_2 <...< p_{\frac{\lg n }{\lg \lg n}}$ are prime numbers such that for every $k$, $p_{2k-1}, p_{2k}= \Theta (n^k)$ and also $\frac{p_{2k+1}}{ p_{2k}} > n$. For example, let $(1+ i\epsilon ) n^{\lfloor\frac {i+1}{2}\rfloor} \leq p_i \leq (1+(i+1)\epsilon) n^{\lfloor\frac{i+ 1}{2}\rfloor}$; where $\epsilon$ is a positive number and $\epsilon=o(\frac{\lg \lg n}{\lg n})$.
Label every edge of $G \setminus H $ randomly and independently by one prime numbers form  $p_1 , \ldots, p_{\frac{\lg n }{\lg \lg n}}$. Next multiply each label of $G \setminus H $ by $2^ {1+ \lceil\lg n \rceil}$.
Assume that $\overrightarrow{O}$ is a random orientation for $G$ such that there is an edge $e=uv$ in $G$ with equal sums of incoming edges at both ends of the edge.

Now, we compute the probability of the event $B_e$ that is "$v$ and $u$ have equal sums of incoming edges".
Let $d_i(v)$ be the number of incoming edges incident with $v$  with label $i$.
Therefore $\sum_i id_i(v) = \sum_i id_i(u)$. According to the definition of prime numbers, for every $k$ we have $d_{2k-1}(v) + d_{2k}(v)= d_{2k-1}(u)+ d_{2k}(u)$ and consequently $d_i(v) = d_i(u)$ for every $i \leq\frac{\lg n}{\lg \lg n}$.
If $e$ is an edge of $H$, since the edges of $H$ have proper edge coloring with distinct powers of two, and since the labels used in $G \setminus H $ are multiplied by $2^ {1+ \lceil \lg n \rceil }$, the event $B_e$ cannot be occur. If $e$ is not in $H$, then at least one of the vertices $u$ and $v$ has degree more than or equal to $\lg n$. Assume that the degree of $v$ is more than or equal to $\lg n$ (note that all of edges incident with $v$, appear in $G \setminus H $).
Let $K=\frac{\lg n}{\lg \lg n}$.
\\ \\
$\Pr(B_e)\leq \Pr\Big(\bigwedge_i (d_i(v) = d_i(u))\Big)=\displaystyle\sum \Pr\Big(\bigwedge_i (d_i(v) = d_i(u)= a_i)\Big).$\\ \\
The above summation is over all vectors $(a_1, \ldots, a_K)$ with $\displaystyle\sum_i a_i \leq d(v)$. \\ \\
\begin{align*}
 \Pr(B_e) &  \leq \displaystyle\sum \Bigg(  \displaystyle\sum_i \Pr \Big( (\overrightarrow{uv} \text{ has label }p_i) \wedge (d_i(u)= d_i(v)+1) \bigwedge_{j, j \neq i} (d_j(u)=d_j(v)) \Big)\\
&\text{ }\text{ }\text{ }\text{ }\text{ }\text{ }\text{ }\text{ }\text{ }\text{ }\text{ }\text{ }+\displaystyle\sum_i \Pr \Big( (\overrightarrow{vu} \text{ has label }p_i) \wedge (d_i(u)= d_i(v)-1) \bigwedge_{j, j \neq i} (d_j(u)=d_j(v)) \Big)\Bigg)\\
&= \displaystyle\sum \Big(\frac{1}{2K}\displaystyle\sum_{a_j \atop j \leq K} \sum_i \Pr(d_j(u)=b_j)\times\Pr(d_j(v)=a_j)+ \frac{1}{2K} \displaystyle\sum_{a_j \atop j \leq K} \displaystyle\sum_i \Pr(d_j(u)=c_j)\times\Pr(d_j(v)=a_j)\Big).
\end{align*}
\\ \\
In the above formula if $j\neq i$, then $b_j= c_j=a_j$   otherwise $b_i=a_i+1$, $c_j=a_i-1$.
\\ \\
\begin{align*}
 \Pr(B_e)  &\leq \displaystyle\sum (1+o(1)) \Pr \Big(\bigwedge_i (d_i(u) =  a_i) \Big)\times \Pr\Big(\bigwedge_i (d_i(v)= a_i)\Big)\\
& \leq \displaystyle\sum (1+o(1)) \Pr \Big(\bigwedge_i (d_i(u) =  a_i) \Big) \times \displaystyle\max_{ a_i \atop i \leq K} \Pr\Big(\bigwedge_i (d_i(v)= a_i)\Big)\\
& = (1+o(1)) \displaystyle\max_{ a_i \atop i \leq K} \Pr\Big(\bigwedge_i (d_i(v)= a_i)\Big)\times \displaystyle\sum  \Pr \Big(\bigwedge_i (d_i(u) =  a_i) \Big)\\
&= (1+o(1)) \displaystyle\max_{a_i\atop i \leq K} \Pr\Big(\bigwedge_i (d_i(v)= a_i)\Big).
\end{align*}
\\
\\
Let $M := \displaystyle \max_{a_i\atop i \leq K} \Pr\Big(\bigwedge_i (d_i(u)= a_i)\Big)$. We have:
\\ \\
\begin{align*}
M &= \Pr \Big(d(v)= \lg n, d^+(v)= \frac{\lg n}{2}, d_1(v)=...=d_{K}(v)\Big)\\
  &= {\displaystyle{\lg n}\choose{\frac{\lg n}{2},\frac{\lg \lg n}{2}, ...,\frac {\lg \lg n}{2}}}(\frac{1}{2})^{\frac{\lg n}{2}}\Big((\frac{1}{2K})^{\frac{\lg \lg n}{2}}\Big)^{K}\\
  &= (\lg n)!\Big((\lg n /2)!((\frac{\lg \lg n}{2})! )^{K}\Big)^{-1}(\frac {1}{4K})^{\frac{\lg n}{2}}.
\end{align*}
\\
\\
By Stirling's approximation $\sqrt{2\pi n} (\frac{n}{e})^n \leq n! \leq \sqrt{e^2 n} (\frac{n}{e})^n$, we have:

\begin{align*}
\frac{(\lg n)!}{(\frac{\lg n}{2})! ((\frac{\lg \lg n}{2})! )^{\frac{\lg n}{\lg \lg n}}}
& \leq \frac{\sqrt{e^2 \lg n} (\lg n /e)^{\lg n}}{(\sqrt{2\pi \frac{\lg n}{2}}(\frac{\lg n}{2e})^{\frac{\lg n}{2}}).(\sqrt{2 \pi \frac{\lg \lg n}{2}} (\lg \lg n /2e)^{\frac{\lg \lg n}{2}})^ \frac{\lg n}{\lg \lg n}}\\
& =(1+o(1))(\frac{4\lg n}{(\lg \lg n)^{1+0.00001}})^{\frac{\lg n}{2}}.
\end{align*}
\\
\\
\\
Consequently,

\begin{align*}
\Pr(B_e) & \leq(1+o(1))M \\
        & = (1+o(1))\Big(\frac {\lg \lg n}{4 \lg n}\Big)^ {\frac{\lg n}{2}}(1+o(1)) \Big(\frac{4\lg n}{(\lg \lg n)^{1+0.00001}}\Big)^{\frac{\lg n}{2}}\\
        & = o(10^{-\frac{\lg n}{2}})\\
        & = o(n^{- \frac{1}{2}}).
\end{align*}

Since the graph $G\setminus H$ contains $\mathcal{O}(n^2)$ edges, so by  linearity of expectation, with the probability $1- o(1)$ the labeling is proper. This completes the proof.

}\end{aliii}

\begin{aliiii}{
{\bf (i)} First, we show that if $G$ is a graph with the universal labeling game number two then $G$ does not contain an odd cycle.  To the contrary suppose that $G$ is graph with some odd cycles. Let $\mathcal{C}=v_1 v_2 \ldots  v_{2k+1} v_1$ be a smallest odd cycle that $G$ contains. The induced graph on the set of vertices $\{v_1,v_2,\ldots, v_{2k+1}\}$ is an odd cycle, otherwise the graph $G$ contains an odd cycle of a size
smaller than $2k+1$. In this graph, if the first player only uses numbers one and two, then the second player has a wining strategy. For each edge $e=v_iu$, $1 \leq i\leq  2k+1$, $u\notin V(\mathcal{C})$ the second player orients $e$ from $v_i$ to $u$. Also for each edge $e=v_iv_{i+1}$ the second player orients $e$ from $v_i$ to $v_{i+1}$. By this approach, without any attention to the  labels of the vertices, the second player wins the game. So if $G$ contains an odd cycle, then $\overrightarrow{\chi_{u}^{g}} (G)>2 $.

Now, let $G$ be a graph without any odd cycle. The graph $G$ is a bipartite graph. If $\Delta(G)>2$, then the second player has a wining strategy. Let $v\in V(G)$, $d(v)>2$ and $v_1,v_2,v_3\in N(v)$. Consider the following strategy for the second player.
For each edge $\{e=vu| u \notin \{v_1,v_2,v_3\}\}$, the second player orients $e$ from $v$ to $u$, Also for each edge
$\{e=v_iu| u\neq v , 1\leq i  \leq 3\}$, the second player orients $e$ from $v_i$ to $u$. The orientations of three edges $vv_1,vv_2,vv_3$ remain unknown (the orientations of other edges are not important). Note that there is no edge between  $v_1,v_2,v_3$ ($G$ is a bipartite graph).
With no loss of generality suppose that the first player chooses $vv_1$ before $vv_2$ and $vv_3$, also chooses $vv_2$ before $vv_3$. If the first player labels $vv_1$ by 1. The second player orients $vv_1$ from $v $ to $v_1$.  Next if the first player labels $vv_2$ by 1. The second player orients $vv_2$ from $v_2$ to $v $, and orients $vv_3$ from $v$ to $v_3 $. otherwise, the second player orients $vv_2$ from $v$ to $v_2 $, and orients $vv_3$ from $v_3$ to $v $. So the second player wins the game. If the first player labels $vv_1$ by 2, by a similar strategy, the second player wins the game.
Therefore, $G$ does not contain any odd cycle and $\Delta(G)\leq 2$. So every connected component of $G$ is an even cycle or a path. It is easy to see  that in this case the first player has a wining strategy.\\
 \\
{\bf (ii)} Let $G$ be a graph with the  set of edges $\{e_1, \ldots, e_{|E(G)|}\}$. We present a strategy for the first player to win the game in $(G, 2\Delta(G))$. For each $i$, the first player chooses  the edge $e_i$ in the round $i$ and  assume that $f_{i-1}:\{e_j: j \in \mathbb{N}_{i-1}\} \rightarrow \mathbb{N}_{2\Delta(G)}$ is the partial labeling which is produced by the first player until the round $i-1$. Also let $D_{i-1}$ be the orientation for the set of edges $\{e_j: j \in \mathbb{N}_{i-1}\}$  which is produced by the second player until the round $i-1$.
Define:

$$S_{i-1}(v)=\displaystyle \sum_{e_j=\overrightarrow{uv}\in D_{i-1} \atop j \leq i-1 } f(e_j).$$

The first player plays such that in each round like $k$, for every two adjacent vertices  $w$ and $z$, $S_{k}(w)\neq S_{k}(z)$.
Now, we present the strategy of the first player in the round $i$.
The first player chooses  the edge $e_i$ in the round $i$, let $e_i=uv$. The first player labels the edge $e_i$ from $\mathbb{N}_{2\Delta(G)}$ such that for every two adjacent vertices $w$ and $z$ of $G$, $S_{i}(w)\neq S_{i}(z)$.
Note that each edge can produce at most one restriction for the value of $f_i(e_i)$.
Thus, in order to make sure that no  two adjacent vertices $w$ and $z$ of $G$ have $S_{i}(w)= S_{i}(z)$, there are at most $2\Delta(G)-1$ restrictions for the value of $f_i(e_i)$. So the first player can find a proper value for $f_i(e_i)$ from $\mathbb{N}_{2\Delta(G)}$, thus has a wining strategy. This completes the proof.\\ \\
{\bf (iii)}
Let $T$ be a tree. First, we show that  $\overrightarrow{\chi_{u}^{g}} (T) \leq \Delta(T)+1$.
Choose an arbitrary vertex $v$ of $T$, and perform a breadth-first search
algorithm from the vertex $v$. This defines a partition $V_0, V_1, \ldots, V_d$ of the vertices of $T$
where each part $V_i$ contains the vertices of $T$ which are at depth $i$ (at distance exactly $i$
from $v$). Let $e_1, \ldots, e_{|E(T)|}$ be an ordering of the edges according to their distance from the vertex $v$.
Define the functions $f_i$, $D_i$ and $S_{i}$ similar to the part $ii$.
The first player plays such that in each round like $k$, for every edge $e_i=wz$, $i\leq k$, we have $S_{k}(w)\neq S_{k}(z)$.
According to the order of the edges and the first player's strategy there are at most $\Delta(T)$ restrictions for $f_k(e_k)$.
Therefore, the first player can find a proper value for $f_k(e_k)$ from $\mathbb{N}_{\Delta(T)+1}$. \\
It is easy to see that the second player has a wining strategy in $(T,\Delta(T)-1)$ (the wining strategy is similar to part $i$). This completes the proof.
}\end{aliiii}

\bibliographystyle{plain}
\bibliography{POref}

\end{document}